\documentstyle[12pt,twoside, epsf]{article}
\def\picill#1by#2(#3)
{\vbox to #2
{\hrule width #1 height 0pt depth 0pt
\vfill\epsffile{#3}}}

\let \ttorg \tt \def \tt{\ttorg \obeyspaces}

\begin{document}

\date{}

\title{\bf A Self-Linking Invariant of Virtual Knots}

\author{Louis H. Kauffman \\
Department of Mathematics, Statistics and
Computer Science\\
University of Illinois at Chicago\\
851 South Morgan St., Chicago IL 60607-7045, USA\\
kauffman@uic.edu
}

 \maketitle
  
 \thispagestyle{empty}

 \subsection*{\centering Abstract}

{\em This paper introduces a self-linking invariant for virtual knots and links, and relates this invariant
to a state model called the binary bracket, and to a class of coloring problems for knots and links
that include classical coloring problems for cubic graphs.} ($AMS$ Subject Classification Number $57M27$)

\section{Introduction}
In this paper we introduce a new invariant of virtual knots and links that is non-trivial for many virtuals, but is trivial on 
classical knots and links. The invariant will initially be expressed in terms of a relative of the bracket polynomial
\cite{bracket}, and then extracted from this polynomial in terms of its exponents, particularly for the case of knots.  This
analog of the bracket polynomial will be denoted $\{ K \}$ (with curly brackets) and called the {\em binary bracket polynomial}.
See Section 3 for the definition and properties of the binary bracket. The key to the combinatorics of this invariant is an
interpretation of the state sum in terms of $2$-colorings of the associated diagrams.
\bigbreak

We define the invariant $$\Lambda(L) = \{ L \}/\Sigma_{O \in O(L)} A^{ w(L^O) }$$ for an unoriented virtual link $L.$
Here $\{ L \}$ denotes the binary bracket, and $\Sigma_{O \in O(L)} A^{ w(L^O) }$ denotes the sum over all orientations of 
$L$ of the terms $A^{w(L^O)}$ where $w(L^O)$ denotes the writhe of $L$ with the specific orientation $O.$ When $\Lambda(L)$ is not
equal to $1$, the virtual link $L$ is non-trivial and non-classical. See Theorem 2 of Section 3.
\bigbreak

In the case of knots, the invariant we extract will be denoted by $J(K)$ for a virtual knot $K.$ It is defined as follows.
Let $w(K)$ denote the writhe of $K.$ ( This is the sum of the crossing signs for any orientation of the knot $K.$)
A crossing $i$ in a knot $K$ is said to be {\em odd} if one encounters an odd number of classical crossings
in walking along the diagram on one full path that starts at $i$ and returns to $i.$ Let $Odd(K)$ denote the set of odd crossings
in the diagram $K.$ In a classical diagram $K,$ the set $Odd(K)$ is empty. 
The invariant $J(K)$ is equal to {\em the sum of the signs of the
crossings in $Odd(K)$}, and we write
$$J(K) = w(K)|_{Odd(K)}.$$ $J(K)$ is invariant under equivalence of virtual knots, and hence is a self-linking number for
virtuals. If $K^{*}$ denotes the mirror image of $K,$ obtained by switching all the classical crossings of $K,$ then 
$J(K^{*}) = - J(K).$ Thus $J(K)$ can detect the difference between virtual knots and their mirror images when it is non-zero.
Since $J(K)$ is zero on classical knots, it detects non-classicality when it is non-zero.  
\bigbreak

The invariants discussed in this paper are elementary. It is particularly striking that the invariant 
$J(K) = w(K)|_{Odd(K)}$ is only infinitesimally more complicated than the classical writhe, and yet can be used to detect
non-triviality, non-classicality and chirality for infinitely many virtual knots.
\bigbreak

The paper is organized as follows. In Section 2 we review facts and definitions about virtual knot theory. The binary bracket is
introduced in Section 3, proofs of invariance, examples and the definitions and theorems about invariants extracted from the
binary bracket are given here. Section 4 delineates a collection of examples of applications of the invariants $J(K)$ and 
$\Lambda(L).$ These include a persistent virtual tangle, a virtual Whitehead link and a virtual Borrommean rings.
Section 5 introduces a combinatorial generalization of the binary bracket to an $n$-ary bracket (that is well-defined on
diagrams, but not an invariant of virtual links) with associated subtle coloring problems for flat virtual shadow diagrams. These
coloring problems are direct generalizations of the
$2$-colorings associated with the binary bracket. We then explore the existence of uncolorables (There  are uncolorable links
even in the case
$n=2.$) for $n$ greater than or equal to $3.$ We show that for $n=3$ the coloring problem defined here is directly related to the
four color theorem in the form of three colorings of the edges of a cubic graph. The section ends with an explanation of the
translation between these subjects.
\bigbreak

\noindent {\bf Acknowledgement.} Much of this effort was sponsored by the Defense
Advanced Research Projects Agency (DARPA) and Air Force Research Laboratory, Air
Force Materiel Command, USAF, under agreement F30602-01-2-05022. 
The U.S. Government is authorized to reproduce and distribute reprints
for Government purposes notwithstanding any copyright annotations thereon. The
views and conclusions contained herein are those of the author and should not be
interpreted as necessarily representing the official policies or endorsements,
either expressed or implied, of the Defense Advanced Research Projects Agency,
the Air Force Research Laboratory, or the U.S. Government. (Copyright 2004.) 
It gives the author great pleasure to acknowledge support from NSF Grant DMS-0245588,
and to give thanks to the University of Waterloo and the Perimeter Institute in Waterloo, 
Canada for their hospitality during the preparation of this research. The author is pleased to thank
Bruce Richter for a helpful conversation about perfect matching.
\bigbreak

\section{Virtual Knot Theory}
Virtual knot theory is an extension of classical diagrammatic knot theory. In this extension one adds 
a {\em virtual crossing} (See Figure 1) that is neither an over-crossing
nor an under-crossing.  A virtual crossing is represented by two crossing arcs with a small circle
placed around the crossing point.  
\bigbreak

Moves on virtual diagrams generalize the Reidemeister moves for classical knot and link
diagrams.  See Figures 1, 2 and 3.  One can summarize the moves on virtual diagrams by saying that the classical crossings
interact with one another according to the usual Reidemeister moves. One adds the detour moves for consecutive sequences of
virtual crossings and this completes the description of the moves on virtual diagrams. It is a consequence of  the moves in
Figure 2 that an arc going through any consecutive sequence of virtual crossings 
can be moved anywhere in the diagram keeping the endpoints fixed;
the places  where the moved arc crosses the diagram become new virtual crossings. This replacement
is the {\em detour move}, and is illustrated schematically in Figure 3.  Note that the fourth move in Figure 2 is a local detour
move. (The corresponding moves with two classical crossings and one virtual crossing are not allowed.)
\bigbreak

One way to understand the meaning of virtual diagrams is to regard them as representatives for oriented Gauss codes
(Gauss diagrams) \cite{VKT,GPV}.  
Virtual equivalence is the same as the equivalence relation generated on the collection
of oriented Gauss codes modulo an abstract set of Reidemeister moves on the codes.  The abstract Reidemeister moves on oriented Gauss codes correspond
exactly to Reideimester moves on diagrammatic representations of these codes in the plane (with virtual crossings), plus the use
of the moves on virtual crossings (all consequences of the detour move). These extra moves  make the particular choice of virtual
crossings in a planar representation irrelevant. We know \cite{VKT,GPV} that classical knot theory embeds faithfully in virtual knot theory.
That is, if two classical knots are equivalent through moves using virtual crossings, then they are equivalent solely via
standard  Reidemeister moves.
\bigbreak

\begin{center}
$$ \picill4inby3in(ReidemeisterMoves)  $$
{ \bf Figure 1 -- Reidemeister Moves} \end{center}

\begin{center}
$$ \picill4inby3in(VirtualMoves)  $$
{ \bf Figure 2 -- VirtualMoves} \end{center}  

\begin{center}
$$ \picill4inby2in(LocalDetourMove)  $$
{ \bf Figure 3 -- The Detour Move} \end{center}

Virtual knots have a special  diagrammatic theory that makes handling them
very similar to the handling of classical knot diagrams. With this approach, one can generalize many structures in classical knot
theory to the virtual domain, and use the virtual knots to test the limits of classical problems such as  the question whether
the Jones polynomial detects knots.   Counterexamples to this conjecture exist in the
virtual domain. The simplest example is the code $C = o1+u2+o3-u1+o2+u3-,$ (Here ``o" stands for ``over", ``u" for ``under", plus
and minus signs refer to the orientations of the crossings $1,2,3.$) a virtualized trefoil, non-trivial,  but with unit Jones
polynomial. It is an open problem whether any of these counterexamples are equivalent to  classical knots.  
\bigbreak

There is a useful topological interpretation for this virtual theory in terms of embeddings of links
in thickened surfaces. Regard each 
virtual crossing as a shorthand for a detour of one of the arcs in the crossing through a 1-handle
that has been attached to the 2-sphere of the original diagram.  
By interpreting each virtual crossing in this way, we
obtain an embedding of a collection of circles into a thickened surface  $S_{g} \times R$ where $g$ is the 
number of virtual crossings in the original diagram $L$, $S_{g}$ is a compact oriented surface of genus $g$
and $R$ denotes the real line.  We say that two such surface embeddings are
{\em stably equivalent} if one can be obtained from another by isotopy in the thickened surfaces, 
homeomorphisms of the surfaces and the addition or subtraction of empty handles.  Then we have the
\smallbreak
\noindent
{\bf Theorem \cite{VKT,DVK,CS1}.} (See also \cite{KUP}.) {\em Two virtual link diagrams are equivalent if and only if their
correspondent  surface embeddings are stably equivalent.}  
\smallbreak
\noindent
\bigbreak

\section{The Binary Bracket Polynomial}
In this section we define a variant of the bracket polynomial \cite{bracket}, called the {\em binary bracket polynomial}
and denoted by $\{ K \} = \{ K \}(A)$ for any (unoriented) virtual knot or link $K.$ 
\bigbreak

We first describe the binary bracket as a state summation. In this respect, it has almost exactly the same formalism as
the standard bracket polynomial, except that the value of an unlabeled loop is equal to $2,$ and the loops in each state are
colored with the  colors from the set $\{0,1 \}$ in such a way that the colors appearing at a smoothing are always different.
This restricts the possible states to a very small number and causes the invariant to behave differently on virtual links than
it does on classical links.
\bigbreak

Let $K$ be any unoriented (virtual) link diagram. Define an {\em unlabeled state}, $S$, of $K$  to
be a choice of smoothing for each  classical crossing of $K.$ There are two choices for smoothing a given  crossing, and
thus there are $2^{N}$ unlabeled states of a diagram with $N$ classical crossings.  A {\em labeled state} is a state $S$ such
that the labels $0$ (zero) or $1$ (one) have been assigned to each component loop in the state. 
\bigbreak

$$\picill5inby4in(D12) $$
\begin{center} {\bf Figure 4 - Bracket Smoothings} 
\end{center}
\vspace{3mm}

In  a state we designate each smoothing with $A$ or $A^{-1}$ according to the left-right convention 
shown in Figure 4. This designation is called a {\em vertex weight} of the state. {\em We require of a labeled state
that the two labels that occur at a smoothing of a crossing are distinct.} This is indicated by a bold line between the 
arcs of the smoothing as illustrated in Figure 4.  
Labeled states satisfying this condition at the site of every smoothing will be called {\em properly labeled states.}
If $S$ is a properly labeled state, we let $\{ K|S \}$ denote the product of its vertex weights, and we define the 
two-color bracket polynomial by the equation:   
$$\{ K \} \, = \sum_{S} \{ K|S \}.$$
where $S$ runs through the set of properly labeled states of $K.$

It follows from this definition that $\{ K \}$ satisfies the equations
$$ \picill.5inby.5in(BracketEquation) ,$$ 
$$\{ K \amalg  O \} \, = 2 \{ K \},$$
$$\{ O \} \, =2.$$
  
\noindent The first equation expresses the fact that the entire set of states of a given diagram is
the union, with respect to a given crossing, of those states with an $A$-type smoothing and those
 with an $A^{-1}$-type smoothing at that crossing. In the first equation, we indicate that the colors at the smoothing are
different by the dark band placed between the arcs of the smoothing. The second and the third equations are clear from the formula
defining the state summation. 
\bigbreak
 
The {\em  binary bracket polynomial} , $\{ K \} \, = \, \{ K \}(A)$,  assigns to each unoriented (virtual) link diagram $K$ a 
Laurent polynomial in the variable $A.$
\bigbreak

In computing the  binary bracket, one finds the following behaviour under Reidemeister move I: 
  $$\{ \mbox{\large $\gamma$} \} = A\{\smile\} \hspace {.5in}$$ and 
  $$\{ \overline{\mbox{\large $\gamma$}} \} = A^{-1}\{\smile\} \hspace {.5in}$$

\noindent where \mbox{\large $\gamma$}  denotes a curl of positive type as indicated in Figure 5, 
and  $\overline{\mbox{\large $\gamma$}}$ indicates a curl of negative type, as also seen in this
figure. The type of a curl is the sign of the crossing when we orient it locally. Our convention of
signs is also given in Figure 5. Note that the type of a curl  does not depend on the orientation
we choose.  The small arcs on the right hand side of these formulas indicate
the removal of the curl from the corresponding diagram.  
\bigbreak

Here is the diagrammatic proof of the behaviour of the  binary bracket with a curl in the diagram.
$$ \picill5inby1.3in(BracketCurl) $$
Note that the second diagram contributes zero, since it contains a demand that an arc be colored diffently from itself.
The proof for the opposite curl goes the same way.
\bigbreak

We now make a key observation about the structure of the properly colored states. Note that any link diagram $K$, real or 
virtual, has an underlying 4-regular plane graph $Sh(K),$ that we shall call the {\em shadow} of $K.$ 
It follows (see Figure 4)
from the combinatorics of coloring at a crossing that a properly colored state of a diagram $K$ (virtual or classical) is
the same as a coloring of the edges of the shadow $Sh(K)$ with $0$ and $1$ such that if two edges meet at a vertex and are not 
adjacent in the cyclic order at that vertex, then they receive different colors. This means that a proper coloring is obtained by 
walking along the diagram, crossing at each crossing, with a color change at each classical crossing and no change at each
virtual crossing. It is easy to see that exactly two such colorings exist for any shadow of a knot diagram. If we
have the shadow of a link diagram with an even number of virtual crossings between any two link components, then there are $2^{N}$
proper diagram colorings where
$N$ is the  number of link components. 
\bigbreak

View Figure 7 for an illustration of the coloring statements of the last paragraph for two virtual knots $K$ and $E.$
In this figure, we have labeled one of the colored states on shadow diagrams for each knot. The smoothings that correspond
to these states are indicated by segments drawn through the crossings. Note that in these virtual cases, there are crossings where
the oriented smoothing is different from the smoothing indicated by the colored state. This kind of difference makes it possible
for the invariant to detect some virtual knots.
\bigbreak
  
$$ \picill4.5inby2.5in(D13) $$
\begin{center} {\bf Figure 5 - Crossing Signs and Curls} 
\end{center}

\noindent {\bf Theorem 0.} The  binary bracket is invariant under regular isotopy for virtual links, and it can be  normalized to
an invariant of ambient isotopy by the definition  
$$Inv_{K}(A) = A^{-w(K)} \{K \}(A),$$ where we choose an orientation for $K$, and where $w(K)$ is 
the sum of the crossing signs  of the oriented link $K$. $w(K)$ is called the {\em writhe} of $K$. 
The convention for crossing signs is shown in  Figure 5.  
\bigbreak

\noindent {\bf Proof of Theorem 0.} First we prove invariance under the second Reidemeister move. The diagrammatic proof 
is shown below.
$$ \picill5inby4.7in(BracketTwoMove) $$
Note that in the this expansion the initial second and third terms are zero due to demands for colors to be distinct from 
themselves. In the remaining two terms, the first consists of two arcs connected through an intermediate circle. If the top
arc is colored $X$, then the circle is colored $\sim X$  ($\sim 0 = 1, \sim 1 = 0$) and the bottom arc is hence colored 
$\sim \sim X = X.$ Thus this condition for the first diagram is that the top and bottom arcs have the same color. This is
the same as saying that the arcs of the reversed smoothing have the same color. Combined with the statement that the vertical
arcs in the second diagram are of  different colors, the two diagrams taken together encompass all cases for two vertical arcs.
Hence the invariance under the  second Reidemeister move is proved.
\bigbreak

For invariance under the third Reidemeister move, view Figure 6.

$$ \picill4.5inby3in(FlatThreeMove) $$
\begin{center} {\bf Figure 6 - States for Third Reidemeister Move} 
\end{center}

In Figure 6 we illustrate the general pattern of state labels on the shadows of the two sides of the third Reidemeister 
move. The variables $x, y, z$ take the values $0$ or $1.$ Note that given choices of the values of these variables for the top
free ends of any one of the diagrams, the values on the rest of the diagram are determined by the coloring rule (switch as the 
signal goes through a classical crossing). Thus we need compare only one state at a time before and after the Reidemeister move.
Note further that at a crossing the four labels will be two of one color and two of the other, determining the smoothing
corresponding to the state. If we switch all colors at a given crossing, then the smoothing remains the same. Note that 
before and after the Reidemeister move, corresponding crossings indeed have all their colors switched. The vertex 
weights are determined by the smoothing and therefore the product of the vertex weights is the same before and after the
smoothing. This proves the invariance under the classical Reidemeister third move.
\bigbreak

It remains to prove invariance under the moves involving the virtual crossings. This is quite easy and we leave the details to
the reader. Writhe normalization works to give invariance under all moves because the writhe itself is an invariant of regular
isotopy and invariant under moves involving virtual crossings. This completes the proof of Theorem 0.
$\hfill\Box$
\bigbreak

$$ \picill5inby4in(FirstVirtualExamples) $$
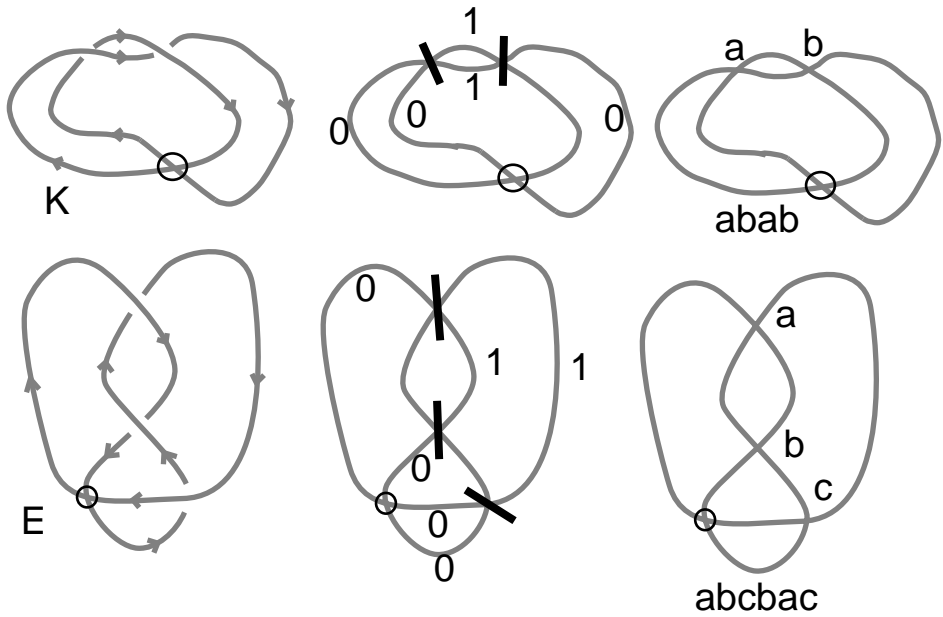
\begin{center} {\bf Figure 7 - Virtual Trefoil $K$ and Virtual Figure Eight $E$} 
\end{center}

\noindent {\bf Remark.} The binary bracket can be viewed as an invariant based on the following solution to the Yang-Baxter
equation.

$$R = \left( \begin{array}{cccc}
0  & 0  & 0  & A \\
0  & A^{-1} & 0 & 0 \\
0 & 0  & A^{-1} & 0 \\
A & 0 & 0 & 0 \\
\end{array} \right).$$

This $4 \times 4$ matrix is viewed as acting upon a tensor product of a two-dimensional space with itself whose basis indices are
$0$ and $1.$ Note that if $A$ is a unit complex number, then $R$ is a unitary matrix. This makes this matrix of interest to 
us in the context of quantum computing as well as topology. See \cite{Dye, BraidGates}.
\bigbreak

\noindent {\bf Theorem 1.} The invariant $Inv$ behaves very simply on classical knot and link diagrams.
\begin{enumerate} 
\item Let $K$ be a classical knot diagram. Then $\{ K \} = 2A^{w(K)}$ where $w(K)$ is the writhe of the diagram
$K$ (Note that for a knot diagram, the writhe is independent of the choice of orientation of that diagram.) Thus
$Inv(K) = 2.$
\item Let $L$ be a classical link diagram. Then $\{ L \} = \Sigma_{O \in O(L)} A^{w(L^O)}$ where $O(L)$ denotes the set of
orientations of $L,$ and $L^{O}$ denotes $L$ with the orientation $O.$ Thus, for a given orientation
$O_{0}$ of L we have $$Inv(L^{O_{0}}) = \Sigma_{O \in O(L)} A^{w(L^O) - w(L^{O_{0}})}.$$
Note that if $L$ has components $\{ L_1, \cdots, L_N \},$ then 
$$w(L^O) - w(L^{O_{0}}) = \Sigma_{i < j} (Lk(L^{O}_{i}, L^{O}_{j}) - Lk(L^{O_{0}}_{i}, L^{O_{0}}_{j})),$$
where $Lk$ denotes the linking number.
\end{enumerate}
\bigbreak

\noindent {\bf Theorem 2.} For virtual diagrams the story is quite different, and $Inv$ can be unequal to $2$ for virtual
knots and unequal to the above writhe or linking number formulas for virtual links.
\begin{enumerate}
\item If a link has an odd number of virtual crossings between two of its components, then there is no proper
coloring of that link diagram. See Figure 8 for an illustration of this in the simplest case of a virtual link $H$ with one
classical crossing and one virtual crossing. The link $H$ has linking number equal to $1/2,$ and linking number alone detects its
linkedness. By convention the value of an empty sum is zero, and hence $\{ H \} = 0,$ whence $Inv(H) = 0.$ Since $Inv(OO) = 4,$
we see that $Inv$ detects the linkedness of $H.$ This case of empty sums is the first example of the use of $Inv$ to detect
virtual links.
\item Call a crossing in a virtual knot diagram $K$ {\em odd} if, in the Gauss code for that diagram there are an odd number of 
appearances of (classical) crossings between the first and the second appearance of $i.$ Let 
$$J(K) = w(K)|_{Odd(K)}$$ where $Odd(K)$ denotes the collection of odd crossings of $K,$ and the restriction of the writhe
to $Odd(K)$, $w(K)|_{Odd(K)},$ means the summation over the signs of the odd crossings in $K.$ Then 
$$Inv(K) = 2A^{-2J(K)}.$$
\item If K is a virtual knot, let $K^{*}$ denote the mirror image for $K$ that is obtained by switching all the crossings of 
the diagram $K.$ Then $$J(K^{*}) = -J(K).$$ Hence, if $J(K)$ is non-zero, then $K$ is inequivalent to its mirror image.
\item If K is a virtual knot and $J(K)$ is non-zero, then $K$ is not equivalent to a classical knot.
\item Let $L$ be a virtual link diagram. 
Let $O(L)$ denote the set of
orientations of $L,$ and $L^{O}$ denote $L$ with orientation $O.$ For a given orientation
$O_{0}$ of $L,$ let  $$\Sigma(L^{O_{0}}) = \Sigma_{O \in O(L)} A^{w(L^O) - w(L^{O_{0}})}.$$
For classical links, $\Sigma(L^{O_{0}}) = Inv(L^{O_{0}}).$ This equality is not always the case for virtual 
links. Nevertheless, $\Sigma(L^{ O_{0} })$ is an invariant of virtual links. When these
two invariants differ, we can conclude that the virtual link is non-trivial and non-classical. The ratio
$$\Lambda(L) = \{ L \}/\Sigma_{O \in O(L)} A^{ w(L^O) } = Inv(L^{O_{0}})/\Sigma(L^{O_{0}}) $$
is an invariant of the unoriented link $L$ that, when not equal to $1,$ determines that the link is non-trivial and 
non-classical.
\end{enumerate}
 
$$ \picill5inby1.5in(H) $$
\begin{center} {\bf Figure 8 - The Link $H$} 
\end{center}

\noindent {\bf Remark.} View Figure 7. The two virtual knots in this figure illustrate the application of Theorem 2. 
In the case of the virtual trefoil $K,$ the Gauss code of the shadow of $K$ is $abab;$ hence both crossings are odd, and we have
$J(K) = 2.$ This proves that $K$ is non-trivial, non-classical and inequivalent to its mirror image. Similarly, the virtual
knot $E$ has shadow code $abcbac$ so that the crossings $a$ and $b$ are odd. Hence $J(E) = 2$ and $E$ is also non-trivial, non-
classical and chiral. Note that for $E,$ the invariant is independent of the type of the even crossing $c.$
\bigbreak

\noindent {\bf Remark.} V. Turaev points out to us \cite{Turaev1} that implicit in the constructions of his paper
on virtual strings \cite{Turaev2} there are interesting generalizations of the invariant $J(K).$  We shall pursue this
topic in another paper.
\bigbreak

\noindent{\bf Proof of Theorem 1.}
\begin{enumerate}
\item To prove the first part, we note that in a classical knot diagram, there is exactly one state and this state has two 
proper colorings. The state can be obtained by choosing one coloring of the diagram and smoothing the crossings accordingly.
Changing all zeros to ones and all ones to zeros gives the other colored state, but does not change the smoothing configuration.
We claim that this smoothing configuration can also be obtained by orienting the diagram and forming an oriented smoothing at each
crossing. (The resulting state is sometimes referred to as the collection of Seifert circles for the diagram.) To see an example,
view Figure 4. The claim follows from the fact that there are an even number of crossings between the first and
second occurrence of any given crossing $i$ in the Gauss code of $K.$ It follows from this that if (say) the color $0$ is the
input color to the crossing $i$, then the color $0$ will also be the output color of the second appearance of the crossing $i.$
The result is that the oriented smoothing of the crossings corresponds to the smoothing designated by the coloring.
Given this claim, we need only point out that the state obtained from the oriented smoothings contributes $A^{w(K)}$ to
the state summation. This follows directly from the definition of the signs of crossings.
\item The proof of this second part requires a generalization of the argument we used in the first part. We need to prove the 
following Lemma.

\noindent {\bf Lemma.} Let $\cal C$ be a collection of Jordan curves in the plane with a set of marked sites (with the structure
of a smoothed crossing). We say that $\cal C$ is {\em properly colored} if  each curve can be assigned the 
color $0$ or the color $1$ such that each site is incident to two distinct colors. Such a proper coloring is possible for $ \cal
C$ if and only if it is possible to orient each Jordan  curve in $\cal C$ such that the orientations at each site are parallel to
one another.

\noindent {\bf Proof of Lemma.} Consider a collection $\cal C$ of oriented Jordan curves in the plane. Each curve has a
well-defined rotation number that is either plus or minus one. By convention, a clockwise oriented circle has rotation
number plus one, while  a counterclockwise oriented circle has rotation number minus one. If $C$ is an oriented Jordan curve,
let $rot(C)$ denote its rotation number. Each curve in $\cal C$ also has a {\em depth} $d(C)$ defined to be 
{\em $d(C) =$ the least number of transverse crossings with curves in $\cal C$ that are needed to draw an arc from the interior
of $C$ to the unbounded region in the plane.} For example, if a curve $C_{1}$ surrounds another curve $C_{2}$ with some pair of 
arcs from the two curves adjacent to one another, then $d(C_2) = 1 + d(C_1).$ In a nest of $n$ circles, the innermost circle has
depth $n$. A curve drawn in the unbounded region has depth $0.$ Now define for each curve $C$ in $\cal C$ the function
$$\lambda(C) = (-1)^{d(C)}rot(C).$$
It is then easy to see that two adjacent curves $C_1$ and $C_2$ in $\cal C$ have parallel orientations if and only if
$\lambda(C_1) \ne \lambda(C_2).$ In Figure 9 we illustrate three curves with locally parallel orientations.
Note that the two concentric curves have the same rotation number, while the two adjacent but not concentric curves have opposite
rotation number.  The Lemma follows from  this observation. 
(We label a curve $C$ with $(1 +\lambda(C))/2$ to change to labels of $0$ and $1$ from labels of $-1$ and $+1.$)
\bigbreak

With this Lemma in hand, we see that every properly colored state of a classical link diagram corresponds to an orientation of
that diagram, and that the evaluation of that state contributes $A$ raised to the writhe of that choice of orientation.
Once the coloring along any given link component is chosen, there is a unique choice of labeling for the rest of the link diagram
to produce a given orientation. The formulas for part $2$ of the Theorem follow directly from these observations.
\end{enumerate}

This completes the proof.$\hfill\Box$
\bigbreak

$$ \picill5inby1.3in(OrientedCircles) $$
\begin{center} {\bf Figure 9 - Nested and Adjacent Oriented Circles} 
\end{center}

\noindent {\bf Proof of Theorem 2.}
Let $K$ be a virtual knot. Just as in the classical case, there are only two labeled states for $K.$ Each state is obtained by 
consecutively labeling the diagram with zeros and ones such that arcs separated by classical crossings are oppositely labeled.
Consider the Gauss code for $Sh(K).$ Without loss of generality, we can assume that the orientation of $K$ is coincident with
the order of the Gauss code.Let $i$ denote one of the classical crossings in $Sh(K).$ We claim that the oriented smoothing at $i$
is identical with the state smoothing at $i$ if and only if the crossing $i$ is even (see the definition of even and odd crossings
given above). To see this claim, view Figure 10. In this Figure we illustrate the case of an even crossing where there are zero
classical crossings between the first and second appearance of $i.$ The local configuration of colors is only changed by changing
the parity of the number of classical crossings between the first and second appearance of $i,$ and we see that in the even case
the state smoothing is coincident with the oriented smoothing.  
\bigbreak

$$ \picill5inby1.3in(EvenCrossing) $$
\begin{center}{\bf Figure 10 - An Even Crossing}
\end{center}
 
We know, therefore, that ${K} = 2A^{a - b}$ where $a$ is the sum of the signs of the even crossings and $b$ is the sum of the
signs of the odd crossings. Note that $w(K) = a + b.$ By definition $J(K) = b$. Hence 
$$Inv(K) = 2A^{a-b -w(K)} = 2A^{a-b-a-b} = 2A^{-2b} = 2A^{-2J(K)}.$$
This completes the proof of the formula  stated in the Theorem. It is clear that changing all the crossings in the knot reverses
the sign of $J(K).$ Since $Inv(K) =2$ for classical knots, we see that $J(K)$ detects non-classicality whenever it is non-zero.
The fifth statement in this Theorem is immediately obvious from the preceding discussion.
This completes the proof of the Theorem. $\hfill\Box$

\section{Examples}
In this section we give a sampling of examples that illustrate the use of the binary bracket polynomial and the associated
self-linking invariant for virtual knots.
\bigbreak

In Figure 11 we show virtual knots $K_1$ and $K_{2}.$ Both knots have underlying flat Gauss code $abcdbdac.$
The code is odd for vertices $a$ and $d$. Thus $J(K_1) = 0$ and $J(K_2 ) =2.$  The invariant $J(K)$ tells us nothing about
$K_1$, but it does tell us that $K_2$ is non-classical and not equivalent to its mirror image. An independent calculation,
that we omit, shows that $K_1$ has unit Jones polynomial, but that it is detected by the two-stand Jones polynomial.
\bigbreak

$$ \picill5inby3in(JPandJQ) $$
\begin{center} {\bf Figure 11 - Two Knots} 
\end{center}

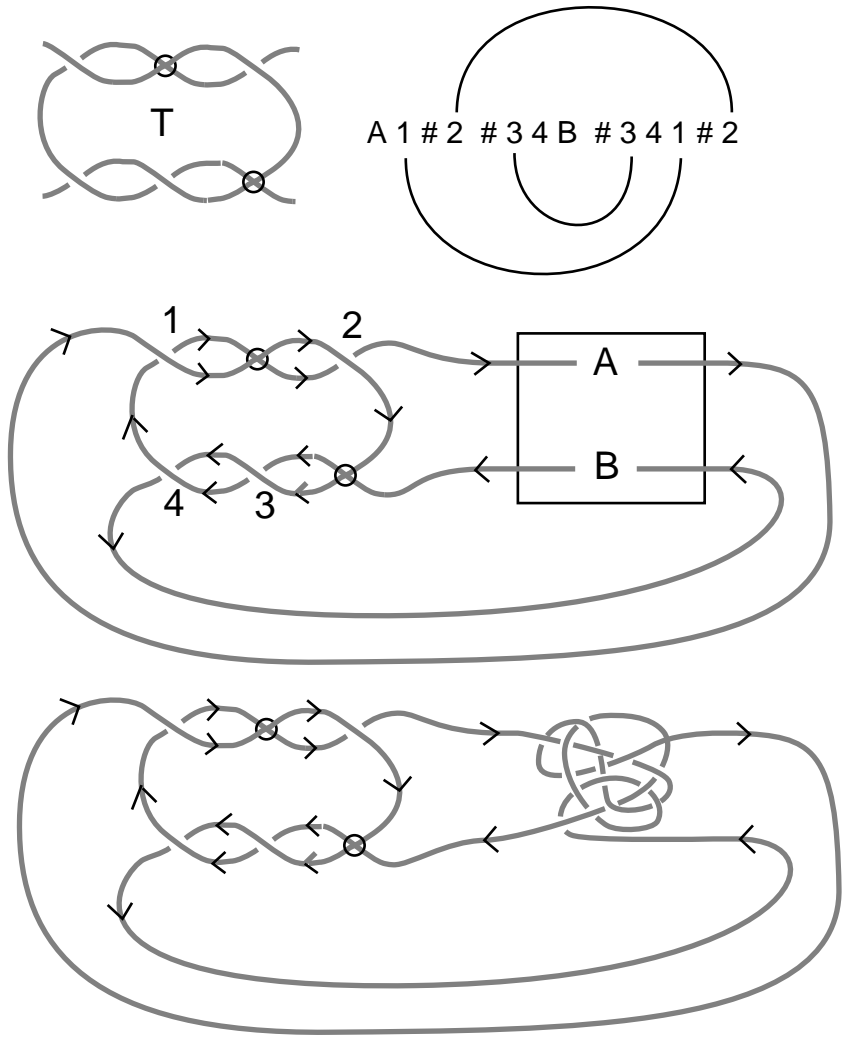
In Figure 12 we illustrate a {\em persistent} virtual tangle $T.$ It follows from the $J$-invariant that whenever this tangle
occurs in a virtual knot diagram $K$ with no other virtual crossings except those in the tangle $T$, then this diagram is
non-trivial, non-classical and inequivalent to its planar mirror image. The proof of this statement is inherent in the Figure. To
see this note that we have indicated a schematic version of  the general code for some diagram in which the tangle sits. The code
has the form $$A1*2*34B*341*2,$$ where the $*$ denotes the occurrence of a virtual crossing in the diagram for the tangle $T,$
and $A$ and $B$ are strings for the remaining part of the Gauss code of $K.$
Since a crossing is odd exactly when its pair of appearances in the Gauss code contains an odd number of virtual crossings, it 
follows that the only odd crossings in the diagram $K$ are $1$, $2$ and $3.$ Hence $J(K) = 3,$ proving the result. 
\bigbreak

$$ \picill5inby6in(PersistentTangle) $$
\begin{center} {\bf Figure 12 - A Persistent Virtual Tangle} 
\end{center}

In Figure 13 we show a ``virtual Whitehead link" $L.$ (The Whitehead link is a classical non-trivial link of two components with
linking number zero.) The link $L$ has $w(L) = -1$ and this is true for each of its four orientations. Hence
$\Sigma_{O \in O(L)} A^{ w(L^O) } = 4A^{-1}.$ The two state diagrams in Figure 13 show that 
$\{L \} = 2(A^{-1} + A^{3}).$ Thus 
$$\Lambda(L) = \{ L \}/\Sigma_{O \in O(L)} A^{ w(L^O) }  = (A^{-1} + A^{3})/2A^{-1} = (1 + A^2)/2.$$ 
Since $\Lambda(L)$ is not equal to $1,$ we conclude that the unoriented link $L$ is not trivial and not classical.
Since $\Lambda(L)(A) \ne \Lambda(L)(A^{-1}),$ we conclude that $L$ is not equivalent to its planar mirror image.
\bigbreak

$$ \picill5inby4in(WhiteheadVirt) $$
\begin{center} {\bf Figure 13 - Virtual Whitehead Link} 
\end{center}

In Figure 14 we show a link $B$ that could be called the ``virtual Borrommean rings." Note that $B$ has writhe zero for each of
its orientations. Thus $\Sigma_{O \in O(B)} A^{ w(B^O) } = 4.$ The states illustrated in the Figure show that
$\{L \} = 2(2 + A^{4} + A^{-4}).$ Thus
$$\Lambda(B) = \{ B \}/\Sigma_{O \in O(B)} A^{ w(B^O) }  = (2 + A^{4} + A^{-4})/2.$$
This shows the the virtual Borrommean rings are non-trivial and non-classical. 
\bigbreak 

$$ \picill5inby6in(VBRings) $$
\begin{center} {\bf Figure 14 - Virtual Borrommean Rings} 
\end{center}

\bigbreak

\section{Colorings and Generalizations}
It is natural to ask what happens in the formalism of the binary bracket if we replace coloring by two colors
with colorings by an arbitrary number of colors. That is, we ask about $\{K\}_n$ where this $n$-ary bracket evalutation
satisfies the equations below.
$$ \picill.5inby.5in(NBracketEquation) ,$$
\smallbreak 
$$\{ K \amalg  O \}_n \, = n \{ K \},$$
$$\{ O \}_n \, = n.$$
Here our conventions are the same as before and $\{K\}_n$ gives a well-defined polynomial  on virtual link diagrams, in the
commuting variables $A$ and $B.$ It appears, however, that unless $n=2$ there is no way to obtain non-trivial 
invariants of (virtual) knots and links from this scheme. Nevertheless, it is of interest to consider the underlying problem of
coloring virtual link diagrams according to the generalization of our rules that is inherent in these equations. To this 
purpose, we define a specialized $n$-ary {\em shadow bracket}, by the following equations.
$$ \picill.5inby.5in(NColorEquation) ,$$ 
\smallbreak
$$[K \amalg  O ]_n \, = n [ K ],$$
$$[ O ]_n \, = n.$$
We call this evaluation of a virtual shadow diagram (note that the crossing is neither over nor under) the shadow bracket to 
emphasize that the crossings in the diagram are flattened. Note that we still have virtual crossings and flat classical crossings.
Coloration at a flat classical crossing follows the rules indicated by the shadow bracket. That is, the rule for coloring is that
{\bf as one crosses a crossing the color must change, and there are exactly two distinct colors at any given flat crossing.} See
Figure 15. At a virtual crossing colors do not change when one crosses the crossing and either one or two colors are present at
the virtual crossing. The equation above for the shadow bracket can be read as tautological. Any coloring at a given crossing
must be in one of the two disjoint possibilities indicated. The value of the shadow bracket on a flat virtual diagram is equal to
the number of colorings of the diagram that are possible under these rules for $n$ colors.
\bigbreak

$$ \picill5inby5in(Rules) $$
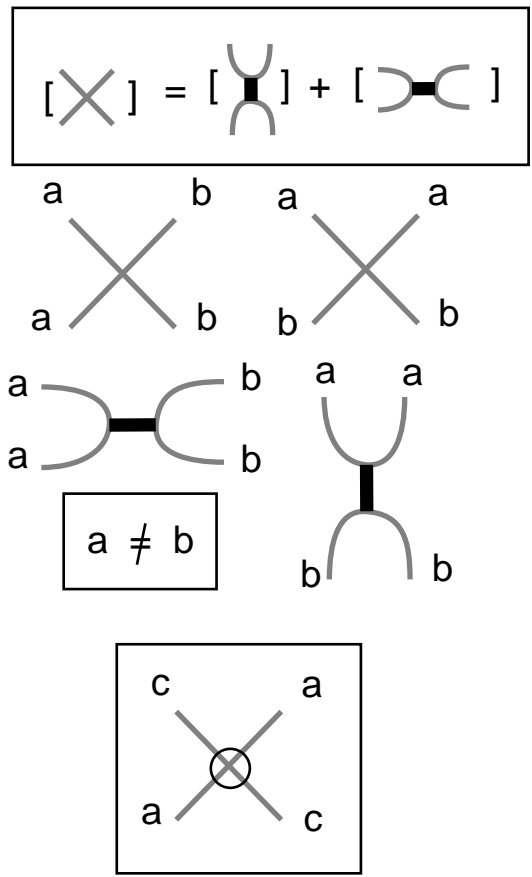
\begin{center} {\bf Figure 15 - Coloring Rules at Flat and Virtual Crossings} 
\end{center}

We would like to know which virtual diagrams are colorable in $n$ colors for $n$ greater than two. When $n$ is equal to two, the
answer is simple, and already used in the virtual knot theory part of this paper. {\em A virtual diagram is colorable with two
colors whenever the number of virtual crossings shared between any two components of the diagram is even.} The situation for 
higher $n$ is much more subtle. First of all consider the diagram in Figure 16. This diagram is the projection of the virtual Hopf
link of Figure 8. It is uncolorable for any $n$, since its structure demands colors that are unequal to themselves.
\bigbreak

$$ \picill5inby2.5in(SimpleUncolorable) $$
\begin{center} {\bf Figure 16 - The Simplest Uncolorable} 
\end{center}

On the other hand, consider the diagram in Figure 17. The reader will have no difficulty in verifying that this diagram can be 
colored in three colors but not in two colors.

$$ \picill5inby2.5in(NeedsThree) $$
\begin{center} {\bf Figure 17 - A Diagram that Needs Three Colors} 
\end{center}

In the top line of  Figure 18 we give an example of a more complex diagram that is uncolorable for any $n.$ Note that an
uncolorable diagram will of necessity have the structure of a flat link diagram that has an odd number of virtual crossings
between some of its components. but it is a subtle matter to characterize the uncolorability.
\bigbreak

One way to begin to understand uncolorables is to look at the expansion of the shadow bracket for one crossing. View Figure 18.
In this figure we illustrate the basic expansion equation for a diagram at one crossing, with the rest of the diagram
concentrated in a tangle box with four external arcs. Uncolorability of $G$ implies that two ways of connecting the arcs on the
tangle box give graphs that, if colorable, force the same color on the two external arcs resulting from the connection.
In the case of the examples shown in Figure 18, it is not hard to see that they satisfy this condition. 
\bigbreak

$$ \picill5inby6in(UncolForm) $$
\begin{center} {\bf Figure 18 - Form of an Uncolorable Diagram} 
\end{center}

The problem of classifying exactly which shadow diagrams are colorable appears to be quite interesting. In fact, it is related 
directly to the classical four color problem \cite{Map,KP}. We now explain this connection.

\subsection{Cubic Graphs, Shadow Diagrams and the Four Color Problem}
A graph is said to be {\em cubic} if there are locally three edges per node. Graphs are allowed to have loops and to have
multiple edges between two nodes. A cubic map $G$ is said to be {\em properly colored} with $n$ colors if the edges of $G$
are colored from the $n$ colors so that all colors incident to any node of $G$ are distinct. It is well known that the following
Theorem is equivalent to the famous Four Color Theorem for maps in the plane.
\bigbreak

\noindent {\bf Theorem} (Equivalent to the Four Color Theorem). Let $G$ be a connected cubic plane graph with no 
isthmus (an {\em isthmus} is an edge whose deletion disconnects the graph). Then $G$ is properly $3$-colorable (as defined above).
\bigbreak

We shall first use this result to give yet another (well-known) equivalent version of the Four Color Theorem (FCT). 
To this end, call a disjoint collection ${\cal E}$ of edges of $G$ that includes all the vertices of $G$ a {\em perfect matching}
of $G.$ Then ${\cal C(E,G)} = G - Interior({\cal E})$ is a collection of {\em cycles} (graphs homeomorphic to the circle, with two
edges incident to each node). We say that ${\cal E}$ is an {\em even} perfect matching of $G$ if every cycle in ${\cal C(E)}$ has
an  even number of edges.
\bigbreak

\noindent {\bf Theorem.} The following statement is equivalent to the Four Color Theorem: Let $G$ be a plane cubic graph with
no isthmus. There there exists an even  perfect matching of $G.$
\bigbreak

\noindent {\bf Proof.} Let $G$ be a cubic plane graph with no isthmus. Suppose that $G$ is properly $3$-colored from the set
$\{ a,b,c \}$. Let ${\cal E}$ denote all edges in $G$ that receive the color $c$. Then, by the definition of proper coloring, the 
edges in ${\cal E}$ are disjoint. By the definition of proper $3$-coloring every node of $G$ is in some edge of ${\cal E}.$
Thus ${\cal E}$ is a perfect matching of $G.$ Since each cycle in ${\cal C(E,G)}$ is two-colored by the the set $\{ a,b \},$
each cycle is even. Hence ${\cal E}$ is an even perfect matching of $G.$
\bigbreak

Conversely, suppose that ${\cal E}$ is an even perfect matching of $G.$ Then we may assign the color $c$ to all the edges of
${\cal E},$ and color the cycles in $\cal{C(E)}$ using $a$ and $b$ (since each cycle is even). The result is a proper
$3$-coloring of  the graph $G.$ This completes the proof of the Theorem. $\hfill\Box$
\bigbreak

\noindent {\bf Remark.} See Figure 19 for an illustration of two perfect matchings of a graph $G.$ One perfect matching is not
even. The other perfect matching is even, and the  corresponding coloring is shown. This Theorem shows that one could
conceivably divide the proving of the FCT into two steps: First prove that every cubic plane isthmus-free graph has a perfect
matching. Then prove that it has an  even  perfect matching. In fact, the existence of a perfect matching is hard, but available
\cite{Petersen}, while the existence of an even perfect matching is really hard! 

$$ \picill5inby5in(EvenOddPartition) $$
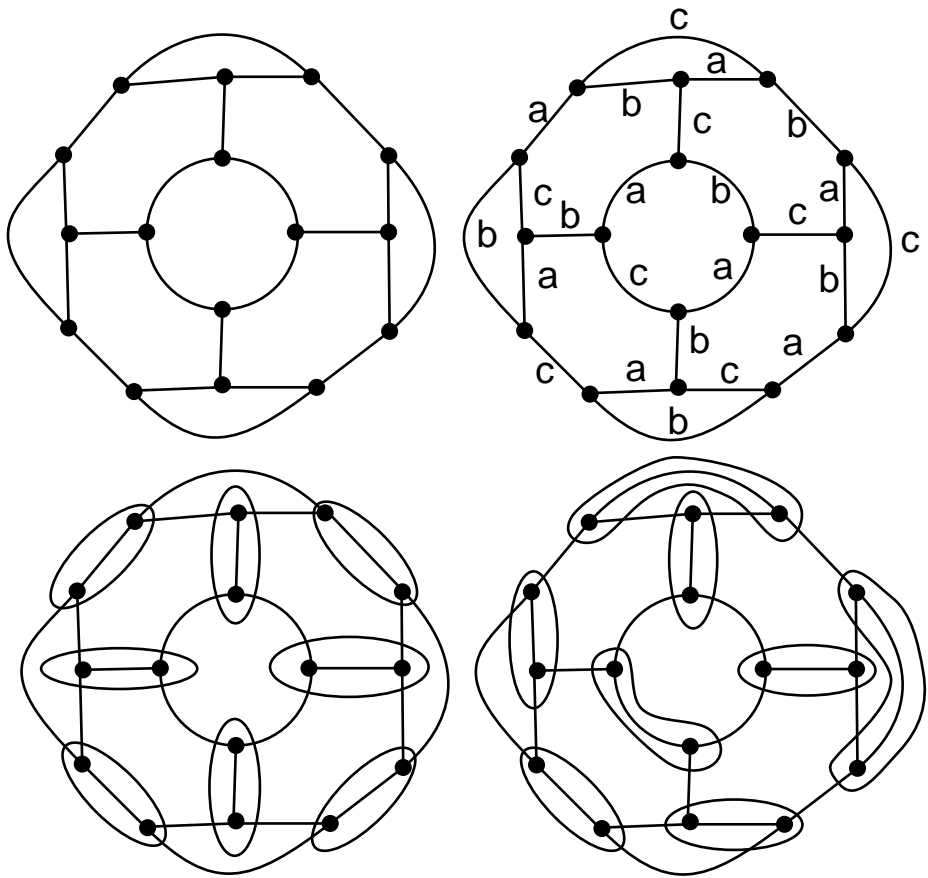
\begin{center} {\bf Figure 19 - Perfect Matchings of a Cubic Plane Graph} 
\end{center}

\noindent {\bf Proposition.} Every cubic graph with no isthmus has a  perfect matching.
\bigbreak

\noindent {\bf Proof.}  See  \cite{Petersen}, Chapter $4.$ $\hfill\Box$
\bigbreak

\noindent {\bf Remark.}  There are graphs that are
uncolorable. Two famous such culprits are indicated in Figure 20. These are examples of graphs with  perfect matchings, but
no even perfect matching. The second example in  Figure 20 is the ``dumbell graph". It is planar, but has an isthmus. The first
example is the Petersen Graph. This graph is non-planar. We have illustrated the Petersen with one  perfect matching that has two
five cycles. No  perfect matching of the Petersen is even. The third  "double dumbell" graph illustrated in Figure 20 has no
perfect matching.
\bigbreak 

$$ \picill5inby4in(DumbPetersen) $$
\begin{center} {\bf Figure 20 - Dumbbell and Petersen} 
\end{center}

We are now in a position to explain the relationship between coloring cubic graphs and coloring virtual shadow diagrams.
Let $G$ be a cubic plane graph without isthmus. We shall say that a graph with no isthmus is {\em bridgeless}.
Let ${\cal E}$ be a perfect matching for
$G.$ Replace each edge in ${\cal E}$ by the combination of flat crossing and virtual crossing shown in Figure 21. Call the
resulting flat virtual diagram $D(G,{\cal E}).$ We have the following Theorem.
\bigbreak

\noindent {\bf Theorem.} Let $G$ be a bridgeless cubic plane graph. Let ${\cal E}$ be a perfect matching for $G.$ Then $G$ can be
properly colored with $3$ colors if and only if $D(G, {\cal E})$ can be properly colored with $3$ colors as a flat virtual
diagram. The Four Color Theorem is equivalent to the statement: There exists a perfect matching for $G$ such that $D(G, {\cal
E})$ can be colored with {\em two} colors. That is, the binary bracket evaluated at $A=1$ does not vanish for $D(G, {\cal E}).$
\bigbreak

$$ \picill5inby4in(Translation) $$
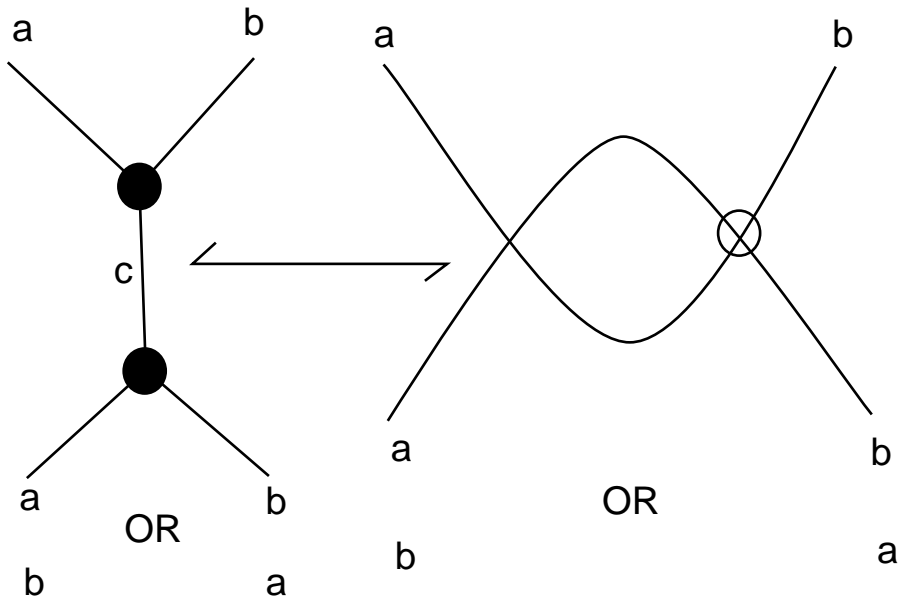
\begin{center} {\bf Figure 21 - Translation between Cubic Graphs and Shadow Virtual Diagrams} 
\end{center}

\noindent {\bf Proof.} As is shown in Figure 21, the coloring conditions for the double $Y$ configuration and for the
replacement shadow diagram are the same when one is coloring at $n=3.$ Note that the edge that is deleted in passing to the
shadow diagram will receive the third color that is distinct from the two colors that appear at the flat crossing in the
shadow diagram. (This shows why this correspondence will not work for $n$ greater than $3.$) Using the perfect matching, one
can replace each edge in
${\cal E}$ with the corresponding shadow diagram configuration. The result is a virtual shadow diagram whose colorings are in  one
to one correspondence with the colorings of the original graph. The rest of the Theorem follows from our discussion of  perfect
matchings and the need for an even perfect matching to satisfy the coloring condition at all vertices of the graph. With an even
perfect matching, the corresponding shadow diagram can be colored with two colors. Hence its binary bracket at $A=1$ does not
vanish.  This completes the proof. 
$\hfill\Box$
\bigbreak

\noindent {\bf Remark.} There is much more to explore in this domain. In Figure 22 we illustrate how the translation
process from cubic graphs to virtual shadow diagrams takes a version of the Petersen graph, with a specific perfect matching to
the uncolorable shadow diagram at the top of Figure 18 (after removal of two redundant virtual crossings). In general, any virtual
shadow diagram can be translated into a cubic graph (with some perfect matching) by placing two canceling virtual crossings next
to any isolated flat crossing in the diagram and then using the combination of flat crossing and virtual crossing to form a
double $Y$ configuration.  The resulting cubic graph may or may not be planar as a result of this operation. For any cubic graph
$G$ with no isthmus, each perfect matching of
$G$ gives rise to a virtual shadow diagram. Thus there is a multiplicity of virtual shadow diagrams corresponding to a given
cubic graph. Note that for $n$ greater than $3$ the colorings for cubic graphs
and the colorings for virtual shadow diagrams are no longer in one-to-one correspondence (since in that case the top and bottom 
ends of the double $Y$ can receive different pairs of colors. 
We shall reserve further comments on this colorful domain for the
next paper.
\bigbreak

$$ \picill5inby4in(PeterToUncol) $$
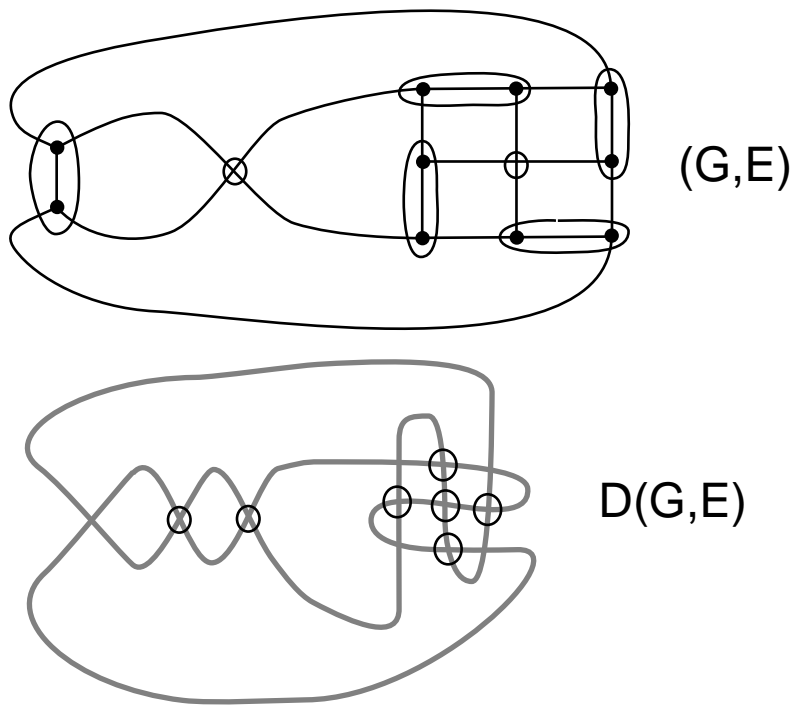
\begin{center} {\bf Figure 22 - Petersen Graph $G$ with perfect matching $E$, and  Virtual Shadow Diagram $D(G,E).$} 
\end{center}

\end{document}